\documentclass[11pt]{article}
\usepackage{amssymb}
\usepackage{multirow}
\usepackage{graphicx}
\usepackage{amsmath,amssymb,verbatim}
\usepackage{amsthm}
\usepackage{epsfig}
\usepackage[latin1]{inputenc}

\usepackage[main=english]{babel}
\newtheorem{Definition}{Definition}[section]

\newtheorem{Theorem}{Theorem}[section]
\newtheorem{Remark}{Remark}[section]
\newcommand{\D}{\mathbb{D}}

\newcommand{\I}{\mathbb{I}}
\newcommand{\R}{\mathbb{R}}
\newcommand{\C}{\mathbb{C}}

\newcommand{\N}{\mathbb{N}}

\begin{document}




\centerline{{\Large{\bf The 1st Level General Fractional
}}}
\vspace{0.2cm}

\centerline{{\Large{\bf 
Derivatives and some of their Properties}}}

\vspace{0.3cm}
\centerline{{\bf Yuri Luchko}}
\vspace{0.3cm}

\centerline{{Department of Mathematics, Physics, and Chemistry}}

\centerline{{Berlin University of Applied Sciences and Technology}}

\centerline{{Luxemburger Str. 10, 13353 Berlin,\ Germany}}

\centerline{{e-mail: luchko@bht-berlin.de}}

\centerline{{telephone number: +49-30-45045295}}

\vspace{0.2cm}

\begin{abstract}
In this paper, we first provide a short summary of the main properties of the so-called general fractional derivatives with the Sonin kernels introduced so far. These are integro-differential operators  defined as compositions of the first order derivative and an integral operator of convolution type. Depending on succession of these operators, the general fractional derivatives of the Riemann-Liouville and of the Caputo types were defined and studied. The main objective of this paper is a construction of the 1st level general fractional derivatives that comprise both the general fractional derivative of the Riemann-Liouville type and the general fractional derivative of the Caputo type. We also provide some of their properties including the 1st and the 2nd fundamental theorems of Fractional Calculus for these derivatives and the suitably defined general fractional integrals.  
\end{abstract}

\vspace{0.2cm}





\section{Introduction}
\label{sec1}

In the framework of the theory of integrals and derivatives  of non-integer order (Fractional Calculus or FC), recently a lot of attention was  attracted by the so called general fractional integrals and derivatives with the Sonin kernels.  These kernels were introduced by Sonin in \cite{Son84} as a natural generalization of the following pair of the power law kernels
\begin{equation}
\label{power}
\kappa(t) = h_\alpha(t),\ \ k(t) = h_{1-\alpha}(t),\ 0<\alpha < 1,
\end{equation}
there $h_\alpha$ is a power law function defined by the formula
\begin{equation}
\label{h}
 h_\alpha(t):=\frac{t^{\alpha-1}}{\Gamma(\alpha)},\ \alpha >0.
\end{equation}

In the papers \cite{Abel1,Abel2}, Abel employed the property 
\begin{equation}
\label{h_con}
(h_\alpha\, *\ h_{1-\alpha})(t) \, = \, \int_0^t \, h_\alpha(t-\tau)\, h_{1-\alpha}(\tau)\, d\tau \, = \, h_1(t)\, = \, \{1\},\ t>0,\ 0<\alpha < 1
\end{equation}
of the kernels $h_\alpha$ and $h_{1-\alpha}$ (by $*$ we denote the Laplace convolution and by $\{1\}$ the function that is identically equal to one for $t>0$) to derive a closed form solution formula for an  integral equation nowadays known  as the Abel integral equation:
\begin{equation}
\label{Abel1}
f(t) = \frac{1}{\Gamma(\alpha)}\, \int_0^t (t-\tau)^{\alpha-1}\phi(\tau)\, d\tau, \ t>0,\ 0<\alpha <1.
\end{equation}
Its solution (derived by Abel formally and in the slightly different notations) has the well-known form:
\begin{equation}
\label{Abel3}
\phi(t) = \frac{d}{dt}\, \frac{1}{\Gamma(1-\alpha)}\, \int_0^t (t-\tau)^{-\alpha}f(\tau)\, d\tau,\ t>0.
\end{equation}

In \cite{Son84}, Sonin used Abel's method to (formally) solve the integral equation 
\begin{equation}
\label{Son1}
f(t) = (\kappa\, *\, \phi)(t)  \, = \, \int_0^t \kappa(t-\tau)\phi(\tau)\, d\tau,\ t>0
\end{equation}
with a kernel $\kappa$  that possesses an associated kernel $k$ satisfying the following condition (compare it to the property \eqref{h_con}): 
\begin{equation}
\label{Son}
(\kappa \, *\, k )(t) = \int_0^t \, \kappa(t-\tau)\, k(\tau)\, d\tau\, = \, h_1(t)\,= \, \{1\},\ t>0. 
\end{equation}
The relation \eqref{Son} is nowadays  referred to as the Sonin condition and the kernels $\kappa$ and $k$ satisfying this condition are called the Sonin kernels. 

Sonin's solution to the equation \eqref{Son1} looks very similar to Abel's solution \eqref{Abel3} with the only difference that the power law kernel $h_{1-\alpha}$ is replaced by the kernel $k$ associated to the kernel $\kappa$:
\begin{equation}
\label{Son2}
\phi(t) = \frac{d}{dt} (k\, *\, f)(t)  \, = \, \frac{d}{dt}\, \int_0^t k(t-\tau)f(\tau)\, d\tau,\ t>0.
\end{equation}

Another important contribution of Sonin to this subject was a construction of a  general class of the Sonin kernels in form of the products of the power law functions and analytical functions. In particular, we mention a famous pair of the Sonin  kernels that was derived in \cite{Son84}:
\begin{equation}
\label{Bess}
\kappa(t) = (\sqrt{t})^{\alpha-1}J_{\alpha-1}(2\sqrt{t}),\
k(t) = (\sqrt{t})^{-\alpha}I_{-\alpha}(2\sqrt{t}),\ 0<\alpha <1,
\end{equation}
where $J_\nu$ and $I_\nu$ are the Bessel and the modified Bessel functions, respectively.

In the modern FC, the right-hand side of the Abel integral equation \eqref{Abel1} is referred to as the Riemann--Liouville fractional integral of the order $\alpha,\ \alpha >0$
\begin{equation}
\label{RLI}
(I^\alpha_{0+}\, f)(t) := (h_\alpha\, * \, f)(t) \, = \, \frac{1}{\Gamma(\alpha)}\, \int_0^t (t-\tau)^{\alpha-1}f(\tau)\, d\tau,\ t>0
\end{equation}
and the right-hand side of the Abel solution formula \eqref{Abel3} is called the Riemann--Liouville fractional derivative of the order $\alpha,\ 0 \le \alpha<1$:
\begin{equation}
\label{RLD}
(D^\alpha_{0+}\, f)(t) := \frac{d}{dt}\, (I^{1-\alpha}_{0+}\, f)(t),\ t>0.
\end{equation}

For $\alpha=0$, the Riemann--Liouville fractional integral is defined as the identity operator (of course, this definition is not arbitrary, but follows from the properties of the operator $I^\alpha_{0+}$ as $\alpha \to 0+$):
\begin{equation}
\label{RLI_0}
(I^0_{0+}\, f)(t) := f(t),\  t>0.
\end{equation}

In its turn, the formula
\begin{equation}
\label{RLD_0}
(D^0_{0+}\, f)(t) \, = \,\frac{d}{dt}\, (I^{1}_{0+}\, f)(t)\, = \, f(t),\ t>0
\end{equation}
is nothing else as the 1st fundamental theorem of calculus. 

In the framework of the Sonin kernels, the above relations are direct consequences from the Sonin condition  \eqref{h_con} extended to the case $\alpha = 0$:
\begin{equation}
\label{h_con_0}
h_0\, *\ h_1  = h_1.
\end{equation}
Of course, the formula \eqref{h_con_0} has to be interpreted in the sense of the generalized functions (the function $h_0$ plays the role of the Dirac $\delta$-function). 

It is worth mentioning that in the theory of the fractional differential equations, another kind of the fractional derivatives is nowadays often employed ($0 \le \alpha < 1$):
\begin{equation}
\label{CD}
( _*D^\alpha_{0+}\, f)(t) :=  (I^{1-\alpha}_{0+}\, f^\prime)(t) \, = \, (D^\alpha_{0+}\, f)(t) - f(0)h_{1-\alpha}(t),\ t>0.
\end{equation}
This operator is usually referred to  as the Caputo fractional derivative of the order $\alpha,\ 0\le \alpha<1$ even if it was introduced and treated by several mathematicians before publications by Caputo and even if Caputo did not in fact consider this form of the operator and worked with it in the Laplace domain.

One more form of the fractional derivatives became recently a subject of active research in FC and its applications. This derivative is defined as follows:  
\begin{equation}
\label{Hil}
(D^{\alpha,\beta}_{0+}\, f)(t)  = (I_{0+}^{\beta(1-\alpha)}\, \frac{d}{dt}\, I_{0+}^{(1-\alpha)(1-\beta)}\, f)(t),\ \ 
0 \le \alpha < 1, \ 0 \le \beta \le 1.
\end{equation}
The operator \eqref{Hil} is usually referred to as the generalized Riemann--Liouville fractional derivative or the Hilfer fractional derivative of order $\alpha$ and type $\beta$ even if it is a particular case of the fractional derivatives introduced by Djrbashian and Nersessian in \cite{DN}, see also \cite{Luc20}. For the properties of this fractional derivative we refer to \cite{Hil00} and subsequent publications \cite{Hil19, Hil, KocLuc19}. 

For $\beta = 0$, this fractional derivative is reduced to the Riemann--Liouville fractional derivative of order $\alpha$ whereas for $\beta =1$ it coincides with the Caputo fractional derivative of order $\alpha$:
\begin{equation}
\label{Hil_RL}
(D^{\alpha,0}_{0+}\, f)(x)  = (I_{0+}^{0}\, \frac{d}{dt}\, I_{0+}^{1-\alpha}\, f)(t)\, = \, \frac{d}{dt}\, (I_{0+}^{1-\alpha}\, f)(t)\, = \,(D^\alpha_{0+}\, f)(t),\ 0 \le \alpha < 1,
\end{equation}
\begin{equation}
\label{Hil_C}
(D^{\alpha,1}_{0+}\, f)(t) \, = \,  (I_{0+}^{1-\alpha}\, \frac{d}{dt}\, I_{0+}^{0}\, f)(t)\, = \, (I_{0+}^{1-\alpha}\, \frac{d}{dt}\, f)(t)\, = \,( _*D^\alpha_{0+}\, f)(t),\ 
0 \le \alpha < 1.
\end{equation}

Thus, treating the fractional derivative $D^{\alpha,\beta}_{0+}$ and the differential equations with this derivative allows avoiding a duplication of work, namely, a separate consideration of the Riemann--Liouville  and the Caputo fractional derivatives and the fractional differential equations with these derivatives.  

In \cite{Koch11}, Kochubei interpreted the operators \eqref{Son1} and \eqref{Son2} with the Sonin kernels $\kappa$ and $k$ as the general fractional integral (GFI) and the general fractional derivative (GFD), respectively: 
\begin{equation}
\label{GFI}
(\I_{(\kappa)}\, f)(t) := (\kappa\, *\, f)(t) \, = \, \int_0^t \kappa(t-\tau)f(\tau)\, d\tau,\ t>0,
\end{equation}
\begin{equation}
\label{FDR-L}
(\D_{(k)}\, f)(t) := \frac{d}{dt} (k\, *\, f)(t) = \frac{d}{dt}\, (\I_{(k)}\, f)(t),\ t>0.
\end{equation}

Kochubei defined also another type of the GFD in the form
\begin{equation}
\label{FDC}
( _*\D_{(k)}\, f) (t) :=  (\D_{(k)}\, f) (t) - f(0)k(t),\ t>0.
\end{equation}

For a function $f$ whose first order derivative is integrable, the operator $ _*\D_{(k)}$ can be represented as follows:
\begin{equation}
\label{FDC_1}
( _*\D_{(k)}\, f) (t) =  (\I_{(k)}\, f^\prime)(t) ,\ t>0.
\end{equation}
In the case of the power law kernel $k(t) = h_{1-\alpha}(t),\ 0\le \alpha <1$, this representation is often used as a definition of the Caputo fractional derivative of the order $\alpha$. 

In \cite{Koch11} and in the subsequent publications devoted to the same topic, the integral operator \eqref{GFI} was called the GFI, whereas the integro-differential operators \eqref{FDR-L} and \eqref{FDC} (or \eqref{FDC_1}) were referred to as the GFDs of the Riemann-Liouville and of the Caputo types, respectively. 

In \cite{Koch11}, mainly the GFD of the Caputo type was treated. In particular, this derivative  with the kernels from a special class $\mathcal{K}$ was shown to be a left inverse operator to the  GFI \eqref{GFI} with the Sonin  kernel $\kappa$ associated to the kernel $k$ of the GFD \eqref{FDC}.  Moreover, in \cite{Koch11}, some important properties of  the solutions to the fractional relaxation equation and to  the Cauchy problem for
the  time-fractional diffusion equation with the GFD of the Caputo type \eqref{FDC} were investigated in detail. 

Recently, a series of papers devoted to the GFI and the  GFDs of the Riemann-Liouville and Caputo types with the Sonin kernels that possess an integrable singularity of power function type at the point zero has been published. In \cite{Luc21a}, these operators were studied on the space of functions that are continuous 
 on the real positive semi-axis and have an integrable singularity of power function type at point zero and its suitable sub-spaces. In \cite{Luc21b}, the  GFDs of arbitrary order of the Riemann-Liouville and  Caputo types were introduced and investigated. These operators extend the definitions \eqref{FDR-L} and \eqref{FDC} that correspond to  the ``generalized order'' of the derivatives between zero and one to the case of any positive real order. An important subclass of the kernels of the GFDs of arbitrary order was suggested in \cite{Tar}. The  Cauchy problems for the fractional differential equations  were considered in \cite{Luc21c} in the case of  the GFDs of the Caputo type and in  \cite{Luc21d,Luc22a} in the case of  the GFDs of the Riemann-Liouville type.   We also mention the papers  \cite{Tar1,Tar2,Tar3}, where the theory presented in \cite{Luc21a,Luc21b,Luc21c,Luc21d,Luc22,Luc22a,Tar} was applied for  formulation of a general  fractional dynamics, a general non-Markovian quantum dynamics, and a general fractional vector calculus, respectively. 

The current situation with the GFDs is very similar to the one with the Riemann--Liouville and the Caputo fractional derivatives. Namely, in the publications devoted to this topic, the case of the GFDs of the Riemann-Liouville type and the case of the GFDs of the Caputo type were treated separately from each other. The main result of this paper is a construction of the GFDs that unify and generalize these two important types of the GFDs. Following \cite{Luc20}, where the case of the power law kernels was considered,  we call these operators the 1st level GFDs. These operators are compositions of two special GFIs and one first order derivative.  This construction can be extended to the compositions of the $n+1$  GFIs and $n$ first order derivatives  following the procedure suggested in \cite{DN} for the case of the Riemann-Liouville integrals and then adopted in \cite{Luc20} in connection to the fundamental theorems of FC formulated for these operators. These extended operators are called the $n$th level GFDs; they will be considered elsewhere.

The rest of the paper is organized as follows: In the second section, we provide a short summary of the basic results already derived in the literature for the GFDs of the Riemann-Liouville and  Caputo types. The main results are presented in the third section, where the 1st level GFDs are introduced and investigated. In particular, we formulate and prove the first and the second fundamental theorems of FC for these operators and a suitably defined GFI. We also present a formula for the projector operator for the  1st level GFDs that comprises the form of the natural initial conditions for the fractional differential equations with these derivatives. Finally, in the fourth section, some discussions and directions for further research are presented. In particular, there we suggest a construction of the $n$th level GFDs that will be discussed elsewhere in detail.

\section{Basic properties of the general fractional integrals and derivatives}
\label{sec2}

The GFI \eqref{GFI}, the GFD of the Riemann-Liouville type \eqref{FDR-L}, and the GFD of the Caputo type \eqref{FDC} with the Sonin kernels $\kappa$ and $k$ are a far reaching generalization of the Riemann-Liouville fractional integral \eqref{RLI} and the  Riemann-Liouville and the Caputo fractional derivatives \eqref{RLD} and \eqref{CD}, respectively, with the power law kernels  $\kappa(t) = h_\alpha(t),\ 0<\alpha<1$ and $k(t)=h_{1-\alpha}(t),\ 0<\alpha < 1$. Because of the Sonin condition \eqref{Son}, the "generalized orders" of the GFI \eqref{GFI} and the GFDs \eqref{FDR-L} and \eqref{FDC} with the Sonin kernels $\kappa$ and $k$ are between zero and one. In \cite{Luc21b}, an extension of the Sonin condition as well as the suitable definitions of the GFI and the GFDs of arbitrary order have been suggested. However, in this paper, we restrict ourselves to the case of the  GFI \eqref{GFI} and the GFDs \eqref{FDR-L} and \eqref{FDC} with the Sonin kernels $\kappa$ and $k$. 

Another important remark is that the general framework of the  GFI \eqref{GFI} and the GFD \eqref{FDR-L} has been suggested by Sonin in \cite{Son84} without any connection to FC. Sonin presented his results in the spirit of that time, i.e., without introducing any spaces of functions and without providing conditions for validity of his derivations and the final formulas. In \cite{Koch11},  for the first time, Kochubei treated these matters mathematically rigorously. In particular, he introduced a very important class $\mathcal{K}$ of the Sonin kernels and provided analysis of the operators \eqref{GFI}, \eqref{FDR-L}, and \eqref{FDC} with the kernels from this class  on the suitable spaces of functions. 

In the recent publications \cite{Luc21a,Luc21b,Luc21c,Luc21d,Luc22,Luc22a}, another important class $\mathcal{L}_1$ of the Sonin kernels has been suggested. Moreover, in these publications, the GFI and the GFDs with the kernels from this class have been investigated on the spaces of functions continuous on the positive real semi-axes that can have a power law singularity at the origin and its subspaces. In this paper, we follow and extend this approach to the GFI and the GFDs with the Sonin kernels.   

\begin{Definition}[\cite{Luc21a}]
\label{d_c}
Let the functions $\kappa$ and $k$ satisfy the Sonin condition \eqref{Son} and the inclusions $\kappa \in C_{-1}(0,+\infty)$ and $k\in C_{-1}(0,+\infty)$, where the space of functions $C_{-1}(0,+\infty)$ is defined as follows:
\begin{equation}
\label{C-1}
C_{-1}(0,+\infty)\, := \, \{f:\ f(t)=t^{p}f_1(t),\ t>0,\ p > -1,\ f_1\in C[0,+\infty)\}.
\end{equation}
The set of pairs $(\kappa,\, k)$ of such kernels is denoted by $\mathcal{L}_{1}$. 
\end{Definition}

It is worth mentioning that the kernels of the most time-fractional derivatives and integrals introduced so far belong to the set $\mathcal{L}_{1}$. In particular, the kernels of the Riemann-Liouville fractional integral and derivative, $h_{\alpha}$ and $h_{1-\alpha}$, respectively, are from $\mathcal{L}_{1}$ if $0<\alpha < 1$.

An important subset of $\mathcal{L}_{1}$ was introduced by Sonin in \cite{Son84}:
\begin{equation}
\label{3-3}
\kappa(t) = h_{\alpha}(t) \cdot \, \kappa_1(t),\ \kappa_1(t)=\sum_{k=0}^{+\infty}\, a_k t^k, \ a_0 \not = 0,\ 0<\alpha <1,
\end{equation}
where $\kappa_1=\kappa_1(t)$ has an infinite convergence radius  and
\begin{equation}
\label{3-4}
k(t) = h_{1-\alpha}(t) \cdot k_1(t),\ k_1(t)=\sum_{k=0}^{+\infty}\, b_k t^k,
\end{equation}
where the coefficients $b_k,\ k\in \N_0$ are uniquely determined by the coefficients $a_k,\ k\in \N_0$ as solutions to the following triangular system of linear equations:
\begin{equation}
\label{3-5}
a_0b_0 = 1,\ \sum_{k=0}^n\Gamma(k+1-\alpha)\Gamma(\alpha+n-k)a_{n-k}b_k = 0,\ n\ge 1.
\end{equation}
In particular the kernels \eqref{Bess} provided in terms of the Bessel and the modified Bessel functions 
are from the class $\mathcal{L}_{1}$. Another example of this type was presented in \cite{Zac08}:
\begin{equation}
\label{3-6}
\kappa(t) =
h_{\alpha,\rho}(t),\ 0<\alpha <1,\ \rho\ge 0,
\end{equation}
\begin{equation}
\label{3-7}
k(t) = h_{1-\alpha,\rho}(t) \, +\, \rho\, \int_0^t h_{1-\alpha,\rho}(\tau)\, d\tau,
\end{equation}
where the function $h_{\alpha,\rho}$ is defined by 
\begin{equation}
\label{2-11}
h_{\alpha,\rho}(t) = \frac{t^{\alpha -1}}{\Gamma(\alpha)}\, e^{-\rho t},\ \ \alpha >0,\ \rho \in \R,\ t>0.
\end{equation}
Finally, let us mention a pair of the Sonin kernels from $\mathcal{L}_{1}$ derived in \cite{Han20}:
\begin{equation}
\label{3-8}
\kappa(t) = h_{1-\beta+\alpha}(t)\, +\, h_{1-\beta}(t),\ 0<\alpha < \beta <1,
\end{equation}
\begin{equation}
\label{3-9}
k(t) = t^{\beta -1}\, E_{\alpha,\beta}(-t^\alpha),
\end{equation}
where $E_{\alpha,\beta}$ stands for the two-parameters Mittag-Leffler function defined by the following convergent series:
\begin{equation}
\label{ML}
E_{\alpha,\beta}(z) = \sum_{k=0}^{+\infty} \frac{z^k}{\Gamma(\alpha\, k + \beta)},\ \alpha >0,\ \beta,z\in \C.
\end{equation}

In the rest of this section, we shortly present some properties of the GFI and the GFDs with the kernels from $\mathcal{L}_{1}$ on the space $C_{-1}(0,+\infty)$ and its suitable subspaces. 

The basic properties of the  GFI \eqref{GFI} on the space $C_{-1}(0,+\infty)$ immediately follow from the known properties of the Laplace convolution and a theorem provided below.

\begin{Theorem}[\cite{LucGor99}]
\label{t1}
The   triple   $\mathcal{R}_{-1} = (C_{-1}(0,+\infty),+,*)$     with  the  usual
addition $+$ and  multiplication $*$ in form of  the  Laplace convolution is a commutative ring without divisors of zero.
\end{Theorem}

In particular, the following relations hold true on the space $C_{-1}(0,+\infty)$  (\cite{Luc21a}):
\begin{equation}
\label{GFI-map}
\I_{(\kappa)}:\, C_{-1}(0,+\infty)\, \rightarrow C_{-1}(0,+\infty),\ \kappa \in \mathcal{L}_{1} \ \mbox{(mapping property)},
\end{equation}
\begin{equation}
\label{GFI-com}
\I_{(\kappa_1)}\, \I_{(\kappa_2)} = \I_{(\kappa_2)}\, \I_{(\kappa_1)},\ \kappa_1,\, \kappa_2 \in \mathcal{L}_{1} \ 
\mbox{(commutativity law)}, 
\end{equation}
\begin{equation}
\label{GFI-index}
\I_{(\kappa_1)}\, \I_{(\kappa_2)} = \I_{(\kappa_1*\kappa_2)},\  \kappa_1,\, \kappa_2 \in \mathcal{L}_{1}\ \mbox{(index law)}.
\end{equation}

According to the axioms of FC discussed in \cite{HL19},  a fractional derivative is an operator left-inverse to the corresponding fractional integral. This property is fulfilled both for the GFD of the Riemann-Liouville type \eqref{FDR-L}  and for the GFD of the Caputo type \eqref{FDC}.

\begin{Theorem}[\cite{Luc21a}]
\label{t3}
Let $\kappa \in \mathcal{L}_{1}$ and $k$ be its associated Sonin kernel. 

Then,  the GFD of the Riemann-Liouville type \eqref{FDR-L} is a left inverse operator to the GFI \eqref{GFI} on the space $C_{-1}(0,+\infty)$: 
\begin{equation}
\label{FTL}
(\D_{(k)}\, \I_{(\kappa)}\, f) (t) = f(t),\ f\in C_{-1}(0,+\infty),\ t>0,
\end{equation}
and the GFD of the Caputo type \eqref{FDC} is a left inverse operator to the  GFI  \eqref{GFI} on the space $C_{-1,(k)}(0,+\infty)$: 
\begin{equation}
\label{FTC}
( _*\D_{(k)}\, \I_{(\kappa)}\, f) (t) = f(t),\ f\in C_{-1,(k)}(0,+\infty),\ t>0,
\end{equation}
where 
\begin{equation}
\label{C_k}
C_{-1,(k)}(0,+\infty) := \{f:\ f(t)=(\I_{(k)}\, \phi)(t),\ \phi\in C_{-1}(0,+\infty)\}.
\end{equation}
\end{Theorem}

The statement formulated in Theorem \ref{t3} is usually referred to as the 1st fundamental theorem of FC for the GFDs (see \cite{Luc20} for a discussion of the 1st and the 2nd fundamental theorems of FC for several different kinds of the fractional derivatives). As we see, it is valid both for the GFD of the Riemann-Liouville type \eqref{FDR-L}  and for the  GFD of the Caputo type \eqref{FDC}. However, the spaces of functions, where these GFDs are left-inverse to the GFI  are very different. The same is of course valid for the conventional Riemann-Liouville and Caputo fractional derivatives. It is worth mentioning that the space $C_{-1,(k)}(0,+\infty)$ defined by \eqref{C_k} can be also characterized as follows (\cite{Luc21a}):
$$
C_{-1,(k)}(0,+\infty) = 
 \{f:\ \I_{(\kappa)} f \in C_{-1}^1(0,+\infty)\, \wedge \, (\I_{(\kappa)}\, f)(0) = 0\},
$$
where the space $C_{-1}^1(0,+\infty)$ is defined by 
\begin{equation}
\label{Cn}
C_{-1}^1(0,+\infty) = \{ f\in C_{-1}(0,+\infty):\,  f^\prime \in C_{-1}(0,+\infty) \}.
\end{equation}
This representation of the space $C_{-1,(k)}(0,+\infty)$ is a generalization of the known property 
$$ 
\left\{ f:\, f = I_{0+}^{1-\alpha}\, \phi,\ \phi \in L_1(0,1)\right\}\, = \, 
\left\{ f:\, I_{0+}^{\alpha}f\in \mbox{AC}([0,\, 1])\, \wedge \, (I_{0+}^{\alpha}f)(0) = 0\right\}
$$
of the spaces of functions employed while treating the Abel integral equation (Theorem 2.3 in \cite{Samko}) 
to the case of the  GFI  with the Sonin kernel $\kappa\in \mathcal{L}_{1}$.

The GFD of the Riemann-Liouville type \eqref{RLD} is also a right inverse operator to the GFI  \eqref{GFI} on the space $C_{-1,(\kappa)}(0,+\infty)$, where $\kappa$ is the Sonin kernel associated to the kernel $k$ of the GFD.  Indeed, in this case, we have the representation $f(t) = (\I_{(\kappa)}\, \phi)(t),\ \phi\in C_{-1}(0,+\infty)$ that leads to the following chain of equations:
$$
(\I_{(\kappa)}\, \D_{(k)}\, f) (t) = (\I_{(\kappa)}\, \frac{d}{dt}\, (k\, *\, f) (t) = (\I_{(\kappa)}\, \frac{d}{dt}\, (k\, *\, (\kappa \, *\, \phi)) (t)=
$$
$$
(\I_{(\kappa)}\, \frac{d}{dt}\, (\{1\} *\, \phi)) (t)= (\I_{(\kappa)}\, \phi) (t)= f(t).
$$
 
However, it is not the case if one considers the  GFD of the Riemann-Liouville type \eqref{RLD}  on its natural domain 
\begin{equation}
\label{C1(k)}
C_{-1,(k)}^{1}(0,+\infty)= \{ f\in C_{-1}(0,+\infty):\, \frac{d}{dt}\, \I_{(k)}\, f \in C_{-1}(0,+\infty) \}.
\end{equation}
The inclusion $C_{-1,(\kappa)}(0,+\infty) \subset C_{-1,(k)}^{1}(0,+\infty)$  follows from Theorem \ref{t3}.
For the 1st order derivative, the space
$C_{-1,(k)}^{1}(0,+\infty)$ corresponds to the space of the continuously differentiable functions, whereas the space $C_{-1,(\kappa)}(0,+\infty)$ consists of all functions that can be represented as integrals of continuous functions. 

Now,  we proceed with the 2nd fundamental theorems of FC for the GFD of the Riemann-Liouville type $\D_{(k)}$ and for the  GFD of the Caputo type $\, _*\D_{(k)}$, respectively. 

\begin{Theorem}[\cite{Luc22}]
\label{tgcTf}
Let $\kappa \in \mathcal{L}_{1}$ and $k$ be its associated Sonin kernel. 

For a function $f\in C_{-1,(k)}^{1}(0,+\infty)$, the formula
\begin{equation}
\label{sFTL}
(\I_{(\kappa)}\, \D_{(k)}\, f) (t) =  f(t) - (k\, *\, f)(0)\kappa(t) = f(t) - (\I_{(k)}\, f)(0)\kappa(t),\ t>0
\end{equation}
holds valid.
\end{Theorem}

In the case of the Riemann-Liouville fractional derivative \eqref{RLD} of  order $\alpha,\ 0 < \alpha < 1$, the kernel function $k$ is the power law function $h_{1-\alpha}$ and its associated kernel $\kappa$ is the function $h_\alpha$. The formula \eqref{sFTL} takes the well-known form (see e.g., \cite{Samko}):
\begin{equation}
\label{RL2nd}
    (I_{0+}^\alpha\, D_{0+}^\alpha \, f)(t)\, =\, f(t) - (I_{0+}^{1-\alpha}\, f)(0)h_\alpha(t),\ t>0.
\end{equation}

\begin{Theorem}[\cite{Luc21a}]
\label{t4}
Let $\kappa \in \mathcal{L}_{1}$ and $k$ be its associated Sonin kernel.

Then,  the relation
\begin{equation}
\label{2FTC}
(\I_{(\kappa)}\, _*\D_{(k)}\, f) (t) = f(t)-f(0),\ t>0 
\end{equation}
holds valid on the space $C_{-1}^1(0,+\infty)$ defined as in \eqref{Cn}. 
\end{Theorem}

For the properties of the space $C_{-1}^1(0,+\infty)$ we refer to \cite{LucGor99}. For the Caputo fractional derivative \eqref{CD} of  order $\alpha,\ 0 < \alpha < 1$, the formula \eqref{2FTC} is well-known (see e.g., \cite{LucGor99}):
\begin{equation}
\label{C2nd}
    (I_{0+}^\alpha\, _*D_{0+}^\alpha \, f)(t)\, =\, f(t) - f(0),\ t>0.
\end{equation}

\begin{Remark}
\label{r1}
The results presented in Theorem \ref{tgcTf} and in Theorem \ref{t4} can be reformulated in terms of the projector operators of the GFDs of the Riemann-Liouville and Caputo types, respectively:
\begin{equation}
\label{proj_RL}
(P_{1}\, f)(t) := f(t) - (\I_{(\kappa)}\, \D_{(k)}\, f) (t) =   (I_{0+}^{1-\alpha}\, f)(0)h_\alpha(t),\ t>0,
\end{equation}
\begin{equation}
\label{proj_C}
(P_{2}\, f)(t) := f(t) - (\I_{(\kappa)}\, _*\D_{(k)}\, f) (t) =   f(0),\ t>0.
\end{equation}
The right-hand sides of the formulas \eqref{proj_RL} and \eqref{proj_C} determine the natural initial conditions that should be posed while dealing with the initial-value problems for the fractional differential equations with the GFDs of the Riemann-Liouville and  Caputo types, respectively. 

We refer to \cite{Luc21d, Luc22a} for the results related to the fractional differential equations with the GFDs of the Riemann-Liouville type  and to \cite{Koch11,KK,Luc21c,LucYam20,Sin18} for analysis of the fractional differential equations with the GFDs of the Caputo type.
\end{Remark}

It is also worth mentioning that the relation \eqref{h_con_0} can be employed for an extension of the definitions of the GFDs of the Riemann-Liouville and Caputo types to the case of the kernel $h_0$ and of the definition of the GFI to the case of the kernels $h_0$ and $h_1$ that are the Sonin kernels in the generalized sense:
\begin{equation}
\label{GFI_0}
(\I_{(h_0)}\, f)(t) := (I^0_{0+}\, f)(t) = f(t),\  t>0,
\end{equation}
\begin{equation}
\label{GFI_1}
(\I_{(h_1)}\, f)(t) := (I^1_{0+}\, f)(t),\  t>0,
\end{equation}
\begin{equation}
\label{FDR-L-0}
(\D_{(h_0)}\, f)(t)  := \frac{d}{dt}(\I_{(h_1)}\, f)(t)\, = \, f(t),\ t>0,
\end{equation}
\begin{equation}
\label{CD-0}
( _*\D_{(h_0)}\, f)(t)  := (\I_{(h_1)}\, f^\prime)(t)\, = \, f(t)-f(0),\ t>0.
\end{equation}

\section{The 1st level general fractional derivatives}
\label{sec3}

In this section, we introduce the 1st level GFDs and investigate their basic properties. We start with defining a suitable generalization of the Sonin kernels.

\begin{Definition}
\label{1l_kernel}
The functions $\kappa,\, k_1,\, k_2 \in C_{-1}(0,+\infty)$ that satisfy the condition 
\begin{equation}
\label{1l_cond}
(\kappa\,*\, k_1\, *\, k_2)(t)\, =\, h_1(t)\, = \, \{1\},\, t>0
\end{equation}
are called the 1st level kernels of the GFDs. 
\end{Definition} 

In what follows, the set of the triples $(\kappa,\, k_1,\, k_2)$ of the 1st level kernels will be denoted by $\mathcal{L}_{1}^{1}$ (the denotation $\mathcal{L}_{n}^{m}$ stands for the set of the kernels of the $m$th level GFDs with the generalized order from the interval $(n-1,\ n)$ and we set $\mathcal{L}_{n}^{0}:= \mathcal{L}_{n}$). Because the Laplace convolution is commutative, the triples of the kernels from $\mathcal{L}_{1}^{1}$ are unordered. However, we agree to associate the first kernel $\kappa$ with the GFI defined by \eqref{GFI} (strictly speaking, this operator cannot be called the GFI anymore because its kernel $\kappa$ is not a Sonin kernel) and the kernels $k_1$ and $k_2$ with a suitably defined GFD of the 1st level. Another important remark is that any of the kernels from the triple $(\kappa,\, k_1,\, k_2)$ is unique as soon as two other kernels are fixed. This follows from the property that the ring   $\mathcal{R}_{-1} = (C_{-1}(0,+\infty),+,*)$    does not have any divisors of zero, see Theorem \ref{t1}. The kernel $\kappa$ from a triple $(\kappa,\, k_1,\, k_2)$ will be called the 1st level kernel associated to the pair of the kernels $(k_1,k_2)$. 

The simplest but an important example of the kernels from $\mathcal{L}_{1}^{1}$ is a triple built by the power law functions
\begin{equation}
\label{power_ex}
\kappa(t) = h_\alpha(t),\  k_1(t) = h_\gamma(t),\  k_2(t)=h_{1-\alpha-\gamma}(t),\ t>0
\end{equation}
under the  conditions
\begin{equation}
\label{power_cond}
0<\alpha<1,\ 0<\gamma < 1-\alpha.
\end{equation}
The inequalities \eqref{power_cond} ensure that all three exponents of the power functions from \eqref{power_ex} are from the interval $(-1,\,0)$ and thus these functions possess the integrable singularities at the origin and belong to the space $C_{-1}(0,+\infty)$. For the functions given by \eqref{power_ex}, the condition \eqref{1l_cond} is satisfied  due to the relation
\begin{equation}
\label{alpha_beta}
(h_\alpha\, *\, h_\beta)(t) = h_{\alpha+\beta}(t),\ t>0,\ \alpha>0,\ \beta >0
\end{equation}
that in its turn follows from the well-known connection between the Euler beta- and gamma-functions:
$$
B(\alpha,\beta) = \int_0^1 t^{\alpha-1}(1-t)^{\beta-1}\, dt = \frac{\Gamma(\alpha)\Gamma(\beta)}{\Gamma(\alpha + \beta)}.
$$

Following the procedure suggested by Sonin in \cite{Son84}, a general class of the 1st level kernels can be constructed in form of the products of the power law functions defined by \eqref{power_ex} and the analytical functions. Indeed, we first fix two out of three kernels from the triple $(\kappa,\, k_1,\, k_2)$ in the form
\begin{equation}
\label{3-3-n}
\kappa(t) = h_{\alpha}(t) \cdot \, f_1(t),\ f_1(t)=\sum_{k=0}^{+\infty}\, a_k t^k, \ a_0 \not = 0,\ 0<\alpha <1,
\end{equation}
\begin{equation}
\label{3-3-n-1}
k_1(t) = h_{\gamma}(t) \cdot \, f_2(t),\ f_2(t)=\sum_{k=0}^{+\infty}\, b_k t^k, \ b_0 \not = 0,\ 0<\gamma<1-\alpha.
\end{equation}
For the function $k_2$, the ansatz 
\begin{equation}
\label{3-3-n-2}
k_2(t) = h_{1-\alpha-\gamma}(t) \cdot \, f_3(t),\ f_3(t)=\sum_{k=0}^{+\infty}\, c_k t^k 
\end{equation}
is employed. The unknown coefficients $c_k$  are determined by substituting the power series \eqref{3-3-n}, \eqref{3-3-n-1}, and \eqref{3-3-n-2} into the condition  \eqref{1l_cond} that leads to an infinite system of the linear equations with a triangular coefficient matrix similar to the one presented in  \eqref{3-5}.

Another possibility for producing the 1st level kernels is by using the Laplace transform technique. Let us assume that the Laplace transforms of the kernels $\kappa$, $k_1$, and $k_2$ exist in the complex half-plane $\Re(p) > c$ with a certain constant $c\in \R$. Applying the Laplace transform to the condition \eqref{1l_cond}  leads to the relation
\begin{equation}
\label{Lap}
\tilde{\kappa}(p)\cdot \tilde{k_1}(p) \cdot \tilde{k_2}(p) \, = \, \frac{1}{p},\ \Re(p) > c
\end{equation}
for the Laplace transforms of the 1st level kernels $\kappa$, $k_1$, and $k_2$. Then the tables of the Laplace transforms  and the inverse Laplace transforms can be used to first find some triples that satisfy the equation \eqref{Lap} and then to return back to the time domain. 

For other techniques for construction of the Sonin kernels that can be also applied in the case of the 1st level kernels we refer to \cite{Sam}. 

We also mention that the formula \eqref{h_con_0} leads to the following useful implication for a triple of the 1st level kernels involving the generalized function $h_0$:
\begin{equation}
\label{1l_cond_0}
h_0\,*\, k_1\, *\, k_2\, =\,  h_1 \, \Rightarrow \, k_1\, *\, k_2\, =\,  h_1
\end{equation}
that means that the remaining kernels $k_1$ and $k_2$ build  a pair of the Sonin kernels from $\mathcal{L}_{1}$ (or $\mathcal{L}_{1}^0$). In this sense, the GFDs of the Riemann-Liouville and Caputo types can be interpreted as the $0$th level GFDs. 

Having defined the 1st level kernels, we now proceed with defining the GFDs with these kernels.

\begin{Definition}
\label{1l_GFD_d}
For a triple $(\kappa,\, k_1,\, k_2) \in \mathcal{L}_{1}^{1}$, a 1st level GFD is defined in the form
\begin{equation}
\label{1l_GFD}
( _{1L}\D_{(k_1,k_2)}\, f)(t)\, :=\, \left(\I_{(k_1)}\, \frac{d}{dt}\, \I_{(k_2)}\, f\right)(t),
\end{equation}
whereas the corresponding GFI with the kernel $\kappa$ is given by the right-hand side of the formula \eqref{GFI}.
\end{Definition} 

The main utility of Definition \ref{1l_GFD_d} is that it covers both the GFD of the Riemann-Liouville type and the GFD of the Caputo type. Indeed, taking into account the representation \eqref{GFI_0} we arrive at the following particular cases of the 1st level GFD \eqref{1l_GFD}:
\vspace{0.1cm}

\noindent
I. For $k_1 = h_0$, the GFD of the Riemann-Liouville type with the kernel $k_2$:
\begin{equation}
\label{1l_RLD}
( _{1L}\D_{(h_0,k_2)}\, f)(t) = \left(\I_{(h_0)}\, \frac{d}{dt}\, \I_{(k_2)}\, f\right)(t) = \left(\frac{d}{dt}\, \I_{(k_2)}\, f\right)(t) = (\D_{(k_2)}\, f)(t).
\end{equation}

\vspace{0.1cm}

\noindent
II. For $k_2 = h_0$, the GFD of the Caputo type with the kernel $k_1$:
\begin{equation}
\label{1l_CD}
( _{1L}\D_{(k_1,h_0)}\, f)(t) = \left(\I_{(k_1)}\, \frac{d}{dt}\, \I_{(h_0)}\, f\right)(t) = \left(\I_{(k_1)} \frac{d}{dt}\,  f\right)(t) = ( _*\D_{(k_1)}\, f)(t).
\end{equation}
\vspace{0.1cm}

Please note that the pairs of the remaining kernels $(\kappa,\, k_2)$ and $(\kappa,\, k_1)$ mentioned in the cases I and II, respectively, are the Sonin kernels from $\mathcal{L}_{1}$, see the implication  
\eqref{1l_cond_0}.

Let us now consider the 1st level GFD \eqref{1l_GFD} in the case of the power law kernels from $\mathcal{L}_{1}^{1}$ defined by \eqref{power_ex}. The GFI with a power law kernel is the conventional Riemann-Liouville fractional integral defined by 
\eqref{RLI} and thus the formula \eqref{1l_GFD} takes the following form:
\begin{equation}
\label{1l_GFD_Hilfer}
( _{1L}\D_{(h_\gamma,h_{1-\alpha-\gamma})}\, f)(t)\, =\, \left(I_{0+}^{\gamma}\, \frac{d}{dt}\, I_{0+}^{1-\alpha-\gamma}\, f\right)(t).
\end{equation}
The operator \eqref{1l_GFD_Hilfer} is the Hilfer fractional derivative defined by \eqref{Hil} in a slightly different parametrization introduced in \cite{Luc20} (the representation \eqref{Hil} is the formula \eqref{1l_GFD_Hilfer} with $\gamma = \beta(1-\alpha)$). In  \cite{Luc20}, the form \eqref{1l_GFD_Hilfer} of the Hilfer fractional derivative was called the 1st level fractional derivative. Because the  operator \eqref{1l_GFD} is a natural generalization of the 1st level fractional derivative \eqref{1l_GFD_Hilfer} to the case of the arbitrary 1st level kernels from $\mathcal{L}_{1}^{1}$, we called it the 1st level GFD, see Definition \ref{1l_GFD_d}.

In the previous section, we provided the 1st and the 2nd fundamental theorems of FC for the GFDs of the Riemann-Liouville and  Caputo types. These results were derived in  separate publications and for  different spaces of functions. In this section, we present a unified approach for handling the 1st level GFDs including its particular cases \eqref{1l_RLD} and \eqref{1l_CD} in form of the GFDs of the Riemann-Liouville and  Caputo types.

\begin{Theorem}
\label{t-FT1}
For the kernels $(\kappa,\, k_1,\, k_2) \in \mathcal{L}_{1}^{1}$,  
the 1st level GFD \eqref{1l_GFD} is a left inverse operator to the GFI \eqref{GFI} on the space $C_{-1,(k_1)}(0,+\infty)$ defined by \eqref{C_k}: 
\begin{equation}
\label{FT1L}
( _{1L}\D_{(k_1,k_2)}\, \I_{(\kappa)}\, f) (t) = f(t),\ f\in C_{-1,(k_1)}(0,+\infty),\ t>0.
\end{equation}
\end{Theorem}

\begin{proof}
A function $f\in C_{-1,(k_1)}(0,+\infty)$ can be represented as
$$
f(t) = (\I_{(k_1)}\, \phi)(t) = (k_1\, * \,  \phi)(t)
$$ 
with a function $\phi \in C_{-1}(0,+\infty)$. 
Because the kernels $(\kappa,\, k_1,\, k_2) \in \mathcal{L}_{1}^{1}$ satisfy the condition \eqref{1l_cond}, we arrive at the following chain of the equations
$$
( _{1L}\D_{(k_1,k_2)}\, \I_{(\kappa)}\, f) (t) = \left(\I_{(k_1)}\, \frac{d}{dt}\, \I_{(k_2)}\, \I_{(\kappa)} \, f \right)(t) =  \left(\I_{(k_1)}\, \frac{d}{dt}\, \I_{(k_2)}\, \I_{(\kappa)} \, \I_{(k_1)}\, \phi \right)(t) =
$$
$$
\left(\I_{(k_1)}\, \frac{d}{dt}\, (k_2\, *\, \kappa\, *\, k_1\,*\,  \phi)(t) \right)(t) = 
\left(\I_{(k_1)}\, \frac{d}{dt}\, (\{ 1\}\, * \, \phi)(t) \right)(t) = (\I_{(k_1)}\, \phi)(t) = f(t)
$$
that proves the formula \eqref{FT1L}. 
\end{proof}

For $k_1 = h_0$ and $k_2 = h_0$, the 1st level GFD is reduced to the GFD of the Riemann-Liouville type with the kernel $k_2$  and to the GFD of the Caputo type with the kernel $k_1$, respectively. Thus, Theorem \ref{t-FT1} covers the results presented in Theorem \ref{t3} including the descriptions of the spaces of functions used in its formulation.

In the next theorem (2nd fundamental theorem of FC for the 1st level GFD), we present a formula for a composition of the GFI with the kernel $\kappa$ and the 1st level GFD on the space of functions $C_{-1,(k_2)}^{1}(0,+\infty)$ defined as in \eqref{C1(k)}. 

\begin{Theorem}
For the kernels $(\kappa,\, k_1,\, k_2) \in \mathcal{L}_{1}^{1}$, 
the formula
\begin{equation}
\label{FT2L}
(\I_{(\kappa)}\,  _{1L}\D_{(k_1,k_2)}\, f) (t) = f(t)\, - \, (\I_{(k_2)}\, f)(0)\, (k_1\, *\ \kappa)(t),\ f\in C_{-1,(k_2)}^{1}(0,+\infty),\ t>0
\end{equation}
holds valid.
\end{Theorem}

\begin{proof} 
For a function $f\in C_{-1,(k_2)}^{1}(0,+\infty)$, the 1st level GFD does exist and is from the space $C_{-1,(k_1)}(0,+\infty)$. Let us now introduce an auxiliary function that we denote by $\phi$: 
 \begin{equation}
\label{p1}
\phi(t) := (\I_{(\kappa)}\,  _{1L}\D_{(k_1,k_2)}\, f) (t). 
\end{equation}
Because of the inclusion $ _{1L}\D_{(k_1,k_2)}\, f \in C_{-1,(k_1)}(0,+\infty)$, Theorem \ref{t-FT1} leads to the relation
 \begin{equation}
\label{p2}
( _{1L}\D_{(k_1,k_2)} \, \phi)(t) = ( _{1L}\D_{(k_1,k_2)}\, \I_{(\kappa)}\,  _{1L}\D_{(k_1,k_2)}\, f) (t) =  ( _{1L}\D_{(k_1,k_2)} \, f)(t)
\end{equation}
that implicates  that the function $\phi - f$ belongs to the kernel of the operator $ _{1L}\D_{(k_1,k_2)}$:
 \begin{equation}
\label{p3}
( _{1L}\D_{(k_1,k_2)} \, (\phi-f))(t) = 0,\ t>0.
\end{equation}
Let us now determine the null-space of the 1st level GFD $ _{1L}\D_{(k_1,k_2)}$:
$$
( _{1L}\D_{(k_1,k_2)} \, f)(t) = 0 \, \Leftrightarrow \, \frac{d}{dt}\, \I_{(k_2)}\, f = 0 
\, \Leftrightarrow \, (\I_{(k_2)}\, f)(t) = (k_2\, *\, f)(t)= C \, \Leftrightarrow \, 
$$
$$
(\kappa\, *\, k_1\, *\, k_2\, *\, f)(t)= (\{C\}\, *\, \kappa\, *\, k_1)(t) \, \Leftrightarrow \, 
(\{1\}\, *\, f)(t)= C(\{1\}\, *\, \kappa\, *\, k_1)(t) \, \Leftrightarrow \, 
$$
$$
f(t)= C(\kappa\, *\, k_1)(t).  
$$
Thus, we arrive at the representation
 \begin{equation}
\label{p4}
\mbox{ker}\, _{1L}\D_{(k_1,k_2)} = \{ C(\kappa\, *\, k_1)(t),\ C\in \R \}.
\end{equation}
Combining it with the relation \eqref{p3}, we get the formula
 \begin{equation}
\label{p5}
\phi(t) = f(t) + C(\kappa\, *\, k_1)(t),\ t>0.
\end{equation}
To determine the constant $C$, we apply the operator $\I_{(k_2)}$ to both sides of the representation \eqref{p5}:
 \begin{equation}
\label{p6}
(\I_{(k_2)}\, \phi)(t) = (\I_{(k_2)}\, f)(t) + C(k_2\, *\, \kappa\, *\, k_1)(t) = 
(\I_{(k_2)}\, f)(t) + C.
\end{equation}
On the other hand, we have the following chain of equations:
$$
(\I_{(k_2)}\, \phi)(t) = (\I_{(k_2)}\, \I_{(\kappa)}\, \D_{(k_1,k_2)}\, f) (t) =
(\I_{(k_2)}\, \I_{(\kappa)}\, \I_{(k_1)}\, \frac{d}{dt}\, \I_{(k_2)} f) (t) =
$$
$$
(k_2\, *\, \kappa\, *\, k_1\, * (\frac{d}{dt}\, \I_{(k_2)} f))(t) = (\{ 1\} * (\frac{d}{dt}\, \I_{(k_2)} f))(t) =
$$
$$
(I_{0+}^1 \frac{d}{dt}\, \I_{(k_2)} f)(t) = (\I_{(k_2)}\, f)(t) - (\I_{(k_2)}\, f)(0). 
$$
Comparing the last formula with the formula \eqref{p6}, we  determine the constant $C$ in the form
 \begin{equation}
\label{p7}
 C =  - (\I_{(k_2)}\, f)(0).
\end{equation}
The representation \eqref{FT2L} follows now from the formulas \eqref{p1}, \eqref{p5}, and \eqref{p7} that completes the proof of the theorem. 
\end{proof}

In the case,  the kernel $k_1$ is the (generalized) kernel $h_0$, the 1st level GFD coincides with the GFD of the Riemann-Liouville type with the kernel $k_2$ and the formula \eqref{FT2L} is reduced to the following form (see Theorem \ref{tgcTf}):
$$
(\I_{(\kappa)}\, \D_{(k_2)}\, f) (t) =  f(t) - (\I_{(k_2)}\, f)(0)\kappa(t),\ f\in C_{-1,(k_2)}^{1}(0,+\infty),\ t>0. 
$$

For the kernel $k_2 = h_0$, the 1st level GFD is the  GFD of the Caputo type with the kernel $k_1$ and the formula \eqref{FT2L} takes  the form (see Theorem \ref{t4}):
$$
(\I_{(\kappa)}\, _*\D_{(k)}\, f) (t) = f(t)-f(0),\ f\in C_{-1}^1(0,+\infty),\ t>0. 
$$

Another important particular case of the formula \eqref{FT2L} is the one that corresponds to the Hilfer fractional derivative \eqref{1l_GFD_Hilfer}. In this case, the 1st level kernels $(\kappa,\ k_1,\, k_2)$ are given by the power functions \eqref{power_ex}. Thus, we have the relations
$$
(I_{(\kappa)}\, f)(t) = (I_{0+}^\alpha\, f)(t),\ (\I_{(k_2)}\, f)(t) =(I_{1-\alpha-\gamma}\, f)(t),\ (h_\alpha\, *\, h_\gamma)(t) = h_{\alpha+\gamma}(t)
$$
that lead to the following form of the 2nd fundamental theorem for the Hilfer fractional derivative in form \eqref{1l_GFD_Hilfer}:
\begin{equation}
\label{proj_H}
(I_{0+}^{\alpha}\, _{1L}\D_{(h_\gamma,h_{1-\alpha-\gamma})}\, f )(t) = f(t) - (I_{1-\alpha-\gamma}\, f)(0)\, h_{\alpha+\gamma}(t),\ t>0.
\end{equation}
This formula (in a slightly different parametrization) was derived for the first time in \cite{Hil}. 

\begin{Remark}
\label{r2}
The formula \eqref{FT2L} can be represented in terms of the projector operator of the 1st level GFD as follows:
\begin{equation}
\label{proj_1}
(P_{1L}\, f)(t) := f(t) - (\I_{(\kappa)}\,  _{1L}\D_{(k_1,k_2)}\, f) (t) =   (\I_{(k_2)}\, f)(0)\, (k_1\, *\ \kappa)(t),\ t>0.
\end{equation}
The right-hand side of the formula \eqref{proj_1} determines the form of the natural initial conditions that should be posed while dealing with the initial-value problems for the fractional differential equations with the 1st level GFDs (see \cite{Hil} for an operational method for derivation of the closed form formulas for solutions to the initial-value problems for the linear fractional differential equations with the Hilfer fractional derivatives and the initial conditions determined by the right-hand side of the formula \eqref{proj_H}).
\end{Remark}

\section{Discussions and directions for further research}

In this paper,  a construction of the 1st level GFD that contains both the GFD of the Riemann-Liouville type the GFD of the Caputo type as its particular cases is suggested. To this end, a suitable generalization of a pair of the Sonin kernels to the case of a kernel triple is first introduced. An important particular case of these kernels is given by three power law functions such that the sum of their exponents is equal to one. This case corresponds the so called Hilfer fractional derivative that involves the conventional Riemann-Liouville and Caputo fractional derivatives as its particular cases. 

For the 1st level GFDs, two fundamental theorems of FC are formulated and proved. Based on the 2nd fundamental theorem, a closed form formula for the projector operator of the 1st level GFD is deduced. This operator determines the form of the natural initial conditions while dealing with the initial-value problems containing the 1st level GFDs.

As to the topics for further research, we mention an in-depth investigation of the fractional differential equations with the 1st level GFDs in the linear and non-linear cases (see  \cite{Luc21d, Luc22a} for some results related to the fractional differential equations with the GFDs of the Riemann-Liouville type and \cite{Koch11,KK,Luc21c,LucYam20,Sin18} for analysis of the fractional differential equations with the GFDs of the Caputo type). Because the 1st level GFD contains the GFDs of the Riemann-Liouville and of the Caputo types as its particular cases, one can cover both types of the fractional differential equations while considering the fractional differential equations with the 1st level GFDs.

Another important topic would be an extension of the notions of the 1st level kernels and the 1st level GFDs to the case of an arbitrary order (in this paper, we restricted ourselves to the case of the "generalized  order" between zero and one). To this end, one can follow the schema suggested in \cite{Luc21b} for the case of the GFDs of the Riemann-Liouville and of the Caputo types, namely, first to replace the condition \eqref{1l_cond} with a more general condition
$$
(\kappa\,*\, k_1\, *\, k_2)(t)\, =\, h_n(t)\, = \, \frac{t^{n-1}}{(n-1)!},\ n\in \N
$$
and then to define the 1st level GFD with a kernel pair $(k_1,\, k_2)$ in the form (the corresponding  GFI is generated by the kernel $\kappa$) 
$$
( _{1L}\D_{(k_1,k_2)}\, f)(t)\, :=\, \left(\I_{(k_1)}\, \frac{d^n}{dt^n}\, \I_{(k_2)}\, f\right)(t).
$$

Finally, we mention that the notion of the 1st level GFDs can be generalized to the case of the 2nd and even $n$th level GFDs following the procedure presented in \cite{Luc20} for the case of the Riemann-Liouville fractional integrals. Say, the 2nd level GFDs that correspond to the GFI with the kernel $\kappa$ are compositions of three GFIs and two first order derivatives:
$$
( _{2L}\D_{(k_1,k_2,k_3)}\, f)(t)\, :=\, \left(\I_{(k_1)}\, \frac{d}{dt}\,\I_{(k_2)}\, \frac{d}{dt}\, \I_{(k_3)}\, f\right)(t),
$$
where
$$
(\kappa\,*\, k_1\, *\, k_2\, *\, k_3)(t)\, =\, h_1(t)\, = \, \{1\},\, t>0.
$$
The $n$th level GFDs are defined as follows: 
$$
( _{nL}\D_{(k_1,\dots,k_{n+1})}\, f)(t)\, :=\, \left(\I_{(k_1)}\, \frac{d}{dt}\,\dots\, \I_{(k_n)}
\frac{d}{dt}\, \I_{(k_{n+1})}\, f\right)(t),
$$
where
$$
(\kappa\,*\, k_1\, *\, \dots\, *\, k_{n+1})(t)\, =\, h_1(t)\, = \, \{1\},\, t>0.
$$

All these and further related topics will be considered elsewhere.

\vspace{0.3cm}

\noindent
{\sl Conflict of Interest:} the author declares no conflict of interest 

\vspace{0.3cm}

\noindent
{\sl Funding:} no funding

\vspace{0.3cm}

\noindent
{\sl Ethical Conduct:} not applicable

\vspace{0.3cm}

\noindent
{\sl Data Availability Statements:} not applicable

\end{document}